\documentclass[12pt]{article}

\usepackage{latexsym}

\usepackage{amsmath,amssymb,amsfonts,amscd}

\newcommand\diag{{\rm diag}\;}
\newcommand\Null{{\rm Null}\;}
\newcommand\hmax{{h^*}}
\newcommand\hmin{{h_0}}
\newcommand\bfA{{\bf A}}
\newcommand\bfB{{\bf B}}
\newcommand\bfC{{\bf C}}
\newcommand\bfd{{\bf d}}
\newcommand\bfM{{\bf M}}
\newcommand\bfN{{\bf N}}
\newcommand\bfP{{\bf P}}

\newcommand\bfY{{\bf Y}}

\newtheorem{theorem}{Theorem}[section]

\newtheorem{lemma}[theorem]{Lemma}
\newtheorem{algorithm}[theorem]{Algorithm}
\newtheorem{rem}[theorem]{Remark}

\newenvironment{remark}{\begin{rem}\em}{\end{rem}}
\newenvironment{proof}{{\noindent\bf Proof.\ }}{\qed}

\newcommand{\bA}{{\mathbb A}}

\newcommand{\bC}{{\mathbb C}}

\newcommand\qed{{\hspace*{\fill}$\Box$\vskip12pt plus 1pt}}
\newcommand\sE{{\mathcal E}}
\newcommand\sF{{\mathcal F}}
\newcommand\sG{{\mathcal G}}
\newcommand\sH{{\mathcal H}}

\newcommand\sW{{\mathcal W}}

\newcommand\sZ{{\mathcal Z}}

\newcommand\pn[1]{{\mathbb P}^{#1}}

\newcommand\W{W}
\newcommand\hatW{{\widehat\W}}

\newcommand\abar{{\bar a}}
\newcommand\bbar{{\bar b}}

\begin{document}

\title{An intrinsic homotopy for intersecting
       algebraic varieties}

\author{Andrew J. Sommese\thanks{Department of Mathematics,
University of Notre Dame, Notre Dame, IN 46556-4618, U.S.A.; {\em
Email:} sommese@nd.edu {\em URL:} http://www.nd.edu/{\~{}}sommese.
This material is based upon work supported by the National Science
Foundation under Grant No.\ 0105653; and
  the Duncan Chair of the University of
        Notre Dame.}
\and Jan Verschelde\thanks{Department of Mathematics, Statistics,
and Computer Science, University of Illinois at Chicago, 851 South
Morgan (M/C 249), Chicago, IL 60607-7045, U.S.A.; {\em Email:}
jan@math.uic.edu or jan.verschelde@na-net.ornl.gov {\em URL:}
http://www.math.uic.edu/{\~{}}jan. This material is based upon
work supported by the National Science Foundation under Grant No.\
0105739 and Grant No.\ 0134611.} \and Charles W.
Wampler\thanks{General Motors Research and Development, Mail Code
480-106-359, 30500 Mound Road, Warren, MI 48090-9055, U.S.A.; {\em
Email:} Charles.W.Wampler@gm.com.}}

\date{April 3, 2004}

\maketitle

\begin{abstract}
Recently we developed a diagonal homotopy method to compute a
numerical representation of all positive dimensional components in
the intersection of two irreducible algebraic sets.
In this paper, we rewrite this diagonal homotopy in
intrinsic coordinates, which reduces the number of variables,
typically in half.  This has the potential to save a
significant amount of computation, especially in the iterative
solving portion of the homotopy path tracker. Three numerical
experiments all show a speedup of about a factor two.



\end{abstract}

Our goal is to compute the irreducible decomposition of $A\cap
B\subset\bC^k$, where $A$ and $B$ are irreducible algebraic sets.
In particular, suppose that
\begin{itemize}
  \item $A$ is an irreducible component
of the solution set of a polynomial system $f_A(u)=0$ defined on
$\bC^k$, and similarly
  \item $B$ is an irreducible component of the solution set
of a polynomial system $f_B(u)=0$ defined on $\bC^k$.
\end{itemize}
This includes the important special case when $f_A$ and $f_B$ are
the same system, but $A$ and $B$ are distinct irreducible
components.

Casting this problem into the framework of numerical algebraic
geometry, we assume that all components are represented as
\emph{witness sets}.  For an irreducible component $A\subset\bC^k$
of dimension $\dim(A)$ and degree $\deg(A)$,
a witness set consists of a generic $k-\dim(A)$ dimensional linear
subspace $L\subset \bC^k$ and the $\deg(A)$ points of
intersection $A\cap L$.
We assume that at the
outset we are given such sets for $A$ and $B$, and our goal is to
compute witness sets for the irreducible components of $A\cap B$.
The intersection may break into several such components, and the
components may have various dimensions.  Our methods proceed in
two phases: we first find a witness superset guaranteed to contain
witness points for all the components, then we break this set into
its irreducible components. We recently reported on an algorithm
\cite{SVW6}, herein called the \emph{extrinsic}\footnote{The
terminology extrinsic/intrinsic is in analogy with the homotopies
of~\cite{HSS98}.} homotopy method,
for computing a witness superset for $A\cap B$. This can then be
decomposed into irreducible components using the methods in
\cite{SVW5} and its references.

Abstracting away the details, which are discussed more fully in
\S\ref{Sec:ExtrinsicReview}, the extrinsic method consists of a
cascade of homotopies in unknowns $x\in\bC^N$ and path parameter
$t\in[0,1]$, each of the form
\begin{equation}\label{Eq:ExtrinsicHom}
  H(x,t):=\left[
     \begin{array}{c}
     f(x)\\
     t(P x+p)+(1-t)(Q x + q)
     \end{array}
  \right] = 0
 \end{equation}
where $f:\bC^N\to\bC^m$ is a system of polynomial equations, $P,Q$
are $(N-m)\times N$ full-rank matrices,  and $p,q\in\bC^{(N-m)}$
are column vectors.  There is a homotopy of this form for each
dimension where $A\cap B$ could have one or more solution
components.  We know solution values for $x$ at $t=1$ and wish to
track solution paths $x(t)$ implicitly defined
by~(\ref{Eq:ExtrinsicHom}) as $t\to0$ to get $x(0)$.

At any specific value of $t$, this looks like
\begin{equation}\label{Eq:ExtrinsicHx}
  \widehat H(x,t)=\left[
     \begin{array}{c}
     f(x)\\
     R(t) x+r(t)
     \end{array}
  \right] = 0,
 \end{equation}
where $R=tP+(1-t)Q$ and $r=tp+(1-t)q$.  The homotopy is
constructed such that we are assured that $R(t)$ is full rank for
all $t\in[0,1]$.  Thus, the linear subspace of solutions of
$R(t)x+r(t)=0$ can be parameterized by $u\in\bC^m$ in the form
\begin{equation}\label{Eq:LinSoln}
  x(u,t) = R^\perp(t) u + x_p(t),
\end{equation}
where $x_p(t)$ is any particular solution and $R^\perp(t)$ is the
right null space of $R(t)$, that is, $R^\perp$ is a full-rank
$N\times m$ matrix with $RR^\perp=0$.  We may restrict $\widehat
H$ to this linear subspace to obtain
\begin{equation}\label{Eq:IntrinsicHu}
  \widetilde H(u,t):=\widehat H(x(u,t),t) = f(R^\perp(t) u + x_p(t)) = 0,
\end{equation}
where we have dropped the linear equations because by
construction, they are identically zero for all $t$. We refer to
this as the \emph{intrinsic} form of the equations.

The problem with~(\ref{Eq:IntrinsicHu}) is that it requires
computing $R^\perp$ and $x_p$ at each new value of $t$ as we
follow the homotopy paths. Because of this, $\widetilde H(x)$ offers
little, if any, computational advantage over the extrinsic
$\widehat H(x)$.

Although not generally possible, for some $P,Q,p,q$, one can
convert the extrinsic homotopy~(\ref{Eq:ExtrinsicHom}) into an
intrinsic homotopy of the form
\begin{equation}\label{Eq:LinearIntrinsic}
  \widetilde H(u,t) = f( t(Cu+c)+(1-t)(Du+d) ) =0,
\end{equation}
in which the path parameter $t$ appears linearly.  This means that
the linear algebra to compute $C,D\in\bC^{N\times m}$ and
$c,d\in\bC^N$ is done just once at the outset, rather than being
repeated at each value of $t$. This can save a significant amount
of computation and is also simpler to implement.

This paper is organized as follows.  In
\S\ref{Sec:ExtrinsicReview}, we review the extrinsic homotopies
formulated in \cite{SVW6} for intersecting algebraic varieties,
and in \S\ref{Sec:IntrinsicStart} and
\S\ref{Sec:IntrinsicCascade}, we show how to convert these to the
linear intrinsic form.  A comparison of the numerical behavior of
the extrinsic homotopies and intrinsic homotopies is presented in
\S \ref{Sec:comparison}.

\tableofcontents

\section{Extrinsic Diagonal Homotopies}\label{Sec:ExtrinsicReview}

Let $A \subset \bC^k$ and $B \subset \bC^k$ be as in the opening paragraph,
having dimensions $a$ and $b$ respectively. We have bounds on the
dimension of components of $A\cap B$ as follows. After renaming if
necessary, we may assume $a\ge b$.  The largest possible dimension
of $A\cap B$ is therefore $b$, which happens if and only if $B$ is
contained in $A$.  We can check this possibility using a homotopy
membership~\cite{SVW2} test to see if a generic point of $B$ is in $A$.
If so, we have $A\cap B = B$ and no further computation is needed.
Otherwise, we know that the largest possible dimension of $A\cap
B$ is $b-1$. On the other hand, because the codimension of $A\cap
B$ is at most the sum of the codimensions of the $A$ and $B$, the
smallest possible dimension of any component of $A\cap B$ is
$\max(a+b-k,0)$. For a particular problem, one might have
available some tighter bounds on $\dim(A\cap B)$, and if so, one
can take advantage of that knowledge in the algorithm to follow.
Accordingly, we introduce the symbols $\hmax$ and $h_0$ as
follows:
\begin{align}
    b \ge{} &\hmax > \dim(A\cap B), \label{Eq:defhstar}\\
    \max(a+b-k,0) \le{} &h_0 \le \min(\dim(\mbox{any component of $A\cap
    B$})).\label{Eq:defh0}
\end{align}
Unless we have other knowledge, we use the defaults $\hmax=b$ and
$h_0=\max(a+b-k,0)$.

Instead of working directly in $\bC^k$, we find the intersection
$A\cap B$ by casting the problem into $(u,v)\in\bC^{k+k}$ and
restricting to the diagonal $u-v=0$. More precisely, the product
$X:=A\times B\subset \bC^{k+ k}$ is an affine variety of dimension
$a+b$, i.e., an irreducible affine algebraic set of dimension
$a+b$. The intersection of $A$ and $B$ can be identified, e.g.,
\cite[Ex.\ 13.15]{E} or \cite[pg.\ 122ff]{M}, with $X\cap\Delta$
where $\Delta$ is the diagonal of $\bC^{k+ k}$ defined by the
system
\begin{equation}
  \delta(u,v):=\left[
    \begin{array}{c}
       u_1-v_1\\
       \vdots\\
       u_k-v_k
    \end{array}
  \right]=0
\end{equation}
with $(u,v)$ giving the coordinates of $\bC^{k + k}$.

The initial data consists of witness sets for $A$ and $B$.  That
is, our data for $A$ consists of a generic system $L_A(u)=0$ of
$a$ linear equations and the $\deg(A)$ solutions
$\{\alpha_1,\ldots,\alpha_{\deg(A)}\} \subset \bC^k$ of the system
\begin{equation}\label{alphas}
  \left[
    \begin{array}{c}
       f_A(u)\\
       L_A(u)
    \end{array}
  \right]=0,
\end{equation}
and similarly the data for $B$  consists of a generic system
$L_B(v)=0$ of $b$ linear equations and the $\deg(B)$ solutions
$\{\beta_1,\ldots,\beta_{\deg(B)}\} \subset \bC^m$ of the system
\begin{equation}\label{betas}
 \left[
    \begin{array}{c}
       f_B(v)\\
       L_B(v)
    \end{array}
  \right]=0.
\end{equation}

\begin{remark}
We are not assuming that $A$ and $B$ occur with multiplicity one
in the solution sets of their respective systems $f_A(u)=0$ and
$f_B(v)=0$. If the multiplicity is greater than one, we must use a
singular path tracker~\cite{SVW4}.
\end{remark}

The extrinsic algorithm can be summarized concisely by introducing
a bit of matrix notation. First, let
\begin{equation}\label{Eq:wdef}
  w = \left[
    \begin{array}{c}
       u \\ v
    \end{array}
      \right]
    \in\bC^{2k},
\end{equation}
and introduce a column vector of ``slack'' variables $z\in\bC^k$.
Also, define the $k\times k$ projection matrix
\begin{equation}\label{Eq:Ph}
  \bfP_h = \diag( \underbrace{1,\ldots,1}_h,
     \underbrace{0,\ldots,0}_{k-h} )
\end{equation}
Left multiplication by $\bfP_h$ picks out the first $h$ rows of
its multiplicand and right multiplication picks out the first $h$
columns of its multiplier.  Note also that $\bfP_h^2=\bfP_h$.
Similarly, let $\bfP_{ji}$ be the $k\times k$ matrix
\begin{equation}\label{Eq:Pji}
       \bfP_{ji} = \diag(
     \underbrace{0,\ldots,0}_{j},
     \underbrace{1,\ldots,1}_{i-j},
     \underbrace{0,\ldots,0}_{k-i} ),
\end{equation}
which picks out rows (or columns) $j+1,\ldots,i$. It is useful to
note that $\bfP_j + \bfP_{ji} = \bfP_i$.

The formulation of the homotopy requires several random matrices
as follows.  First, we choose generic matrices
\begin{equation}\label{Eq:RandMats}
  \bfM\in\bC^{(k-a)\times \#(f_A)}, \qquad
  \bfN\in\bC^{(k-b)\times \#(f_B)},
\end{equation}
where $\#(f_A)$ is the number of functions in the system $f_A(x)$
associated to component $A$, and similarly for $\#(f_B)$.  These
are used to define
\begin{equation}\label{Eq:sFw}
 \sF(w) :=
  \left[
    \begin{array}{c}
       \bfM f_A(u)\\
       \bfN f_B(v)
    \end{array}
  \right].
\end{equation}
Note that $A\times B$ is an irreducible component of the solution
set of the system $\sF(w)=0$.  Next, we choose $\bA$ a generic
$(a+b)\times k$ matrix, and let
\begin{equation}\label{Eq:bfA}
  \bfA =
  \left[ \begin{array}{cc} \bA & -\bA \end{array} \right]
  \in\bC^{(a+b)\times 2k}
\end{equation}
so $\bfA w = \bA(u-v)$.  Finally, we choose generic matrices
\begin{equation}\label{Eq:BCDmats}
  \bfB\in\bC^{(a+b)\times k},\qquad \bfC\in\bC^{k\times2k}
  \qquad \bfd\in\bC^{k\times 1}.
\end{equation}
In all these, a matrix with random complex elements will be
generic with probability one.

Since the smallest dimensional nonempty component of $A\cap B$ is
of dimension at least $\max\{0,a+b-k\}$, it follows from
\cite[Lemma (3.1)]{SVW6} that we can find the irreducible
decomposition of $A\cap B$ by finding the irreducible
decomposition of $\bfA w=0$ on $X=A\times B$. For this purpose, we
consider a cascade of homotopies of the form
\begin{equation}\label{Eq:Cascade}
 \sE_{h}(w,z) =
 \left[
    \begin{array}{l}
       \sF(w)\\
       \bfA w+\bfB \bfP_h z\\
       z - \bfP_h(\bfC w + \bfd)
    \end{array}
  \right]=0,
\end{equation}
which is well-defined for any integer $0\le h\le k$. Denoting the
entries of $z$ as $z_1,\ldots,z_k$, note that the last row of this
matrix equation implies that $(z_{h+1},\ldots,z_k)=0$. The method
for generating a witness superset consists of solving
$\sE_\hmax(w,z)=0$ and then descending sequentially down the
cascade to solve $\sE_j(w,z)=0$ for $j=\hmax-1,\ldots,\hmin$.

The rationale behind the cascade is that the linear system
$\bfP_h(\bfC w+\bfd)=0$ is a linear slice that cuts out witness
points for solution components of dimension $h$.  The vector $z$
is a set of slack variables.  A solution point of $\sE_h(w,z)=0$
for which $z=0$ is on the slice and thus gives a witness point.
Solution points with $z\ne0$ are not on the slice, and we call
these ``nonsolutions.'' These become the starting points for the
next step of the cascade. (We state this more formally below,
after giving more details of the algorithm.) For each step down
the cascade, one more slack variable is set to zero and a
corresponding hyperplane is removed from the slice. The recycling
of nonsolutions as starting points for the next step of the
cascade is valid due to the fact that for $j<i$, $\sE_j(w,z)$ is
just $\sE_i(w,z)$ with certain elements of $\bfB$, $\bfC$, and
$\bfd$ set to zero.  This is justified in \cite{SVW6}.

The following steps of the algorithm still need to be described:
\begin{itemize}
  \item how to solve $\sE_\hmax(w,z)=0$,
  \item how to descend the cascade, and
  \item how to reap the witness points from the solutions at each
  level of the cascade.
\end{itemize}
The homotopy to solve $\sE_\hmax(w,z)=0$ is
\begin{equation}\label{Eq:preHomotopy}
  \left[\begin{array}{{c}}
                \sF(w)\\
               (1-t)\left[
    \begin{array}{l}
       \bfA w+\bfB \bfP_\hmax z\\
       z - \bfP_\hmax(\bfC w + \bfd)
    \end{array}
  \right]+t\gamma\left[
                \begin{array}{c}
                 L_A(u)\\
                 L_B(v)\\
                 z
                 \end{array}\right]
               \end{array}
      \right]=0,
\end{equation}
where $\gamma$ is a random complex number. At $t=1$, solution
paths start at the $\deg(A) \times \deg(B)$ nonsingular solutions
 $\{(\alpha_1,\beta_1),\ldots,(\alpha_{\deg(A)},\beta_{\deg(B)})\}
   \subset \bC^{2k}$ obtained by combining the witness points for
$A$ and $B$. At $t=0$, the solution paths terminate at the desired
start solutions for $\sE_\hmax(w,z)=0$.  In \cite{SVW6} we ended
the homotopy at $\sE_b(w,z)=0$, but the argument works equally
well with $\hmax$ in place of $b$.

The homotopy connecting $\sE_i$ to $\sE_j$ for $j<i$ is
\begin{equation}\label{Eq:CascadeHomotopies}
 \sH_{i,j}(\tau,w,z):=
 \left[
    \begin{array}{l}
       \sF(w)\\
       \bfA w + \bfB \bfP_i z\\
       z - (\bfP_j + \tau \bfP_{ji})(\bfC w + \bfd)
    \end{array}
  \right]=0,
\end{equation}
where $\tau$ goes from 1 to 0 along a sufficiently general
$1$-real-dimensional curve. For example, for all but finitely many
$\gamma\in\bC$ of absolute value~$1$, $\tau=r + \gamma r(1-r)$  as
$r$ goes from 1 to 0 on the real interval suffices.  Another
possibility, relevant in what comes below, is
\begin{equation} \label{crazyFormula}
  \tau = t/(t+\gamma(1-t))
\end{equation}
as $t$ goes from 1 to 0 on the real interval.

In the cascade of homotopies from \cite{SVW6} (based on \cite{SV}),
we start out with the finite set $\sG_i$ of nonsingular solutions
of $\sE_i$ with $z_i\not=0$. Tracking these start solutions we end
up with a set of solutions $\sG^\sE_{i,j}$ of $\sE_j$ with $z_h=0$
for $h>j$.  In \cite{SVW6}, $j=i-1$, but the argument there works
immediately for any $j<i$. The key points about the set
$\sG^\sE_{i,j}$ is that
\begin{enumerate}
  \item the set $\sG_j$ equals the set of points in $\sG^\sE_{i,j}$ for
which $z_j\not=0$;
  \item the set of points $\hatW_j \subset \sG^\sE_{i,j}$ for
which $z_h=0$ for all $h\le j$ contains a witness point set $W_j$
for the $j$-dimensional components of the solution set of the
intersection of $A$ and $B$.
\end{enumerate}
We also know that the set of points in $\sG^\sE_{i,j}$ for which
$z_h=0$ for all $h\le j$ equals the set of points in
$\sG^\sE_{i,j}$ for which $z_j=0$. We wish to set up an intrinsic
homotopy such that analogs of the above key facts hold true.

\section{Setting Up Intrinsic Homotopies}

The extrinsic homotopies
of~(\ref{Eq:preHomotopy}) and~(\ref{Eq:CascadeHomotopies}) use the
variables $(w,z)\in\bC^{2k}\times\bC^k$.  Each has $a+b+k$ linear
equations which we wish to eliminate by converting to an intrinsic
homotopy.  The result will be homotopies in intrinsic variables
$y\in\bC^{2k-a-b}$.  Note that $2k-(a+b)$ is the codimension of
$A\times B$ in $\bC^{2k}$.  It is also the sum $\abar+\bbar$ of
the codimension $\abar=k-a$ of $A$ in $\bC^k$ and the codimension
$\bbar=k-b$ of $B$ in $\bC^k$.  Since this quantity appears
frequently in the expressions below, we define
\begin{equation}
  m = 2k-a-b.
\end{equation}
Accordingly, our intrinsic homotopy variables are $y\in\bC^m$.

\subsection{Intrinsic Start Homotopy}\label{Sec:IntrinsicStart}

In this section, we replace the extrinsic start homotopy
of~(\ref{Eq:preHomotopy}) with one having the intrinsic form
of~(\ref{Eq:LinearIntrinsic}).  Fixing a particular solution
\begin{equation}
w_1=\left[
     \begin{array}{c}
           u_p \\
           v_p
     \end{array}\right]
\end{equation}
of
\begin{equation}\label{SysStart}
\left[
                \begin{array}{c}
                 L_A(u)\\
                 L_B(v)
                 \end{array}\right]=0,
\end{equation}
choose a basis $W_1\in\bC^{2k\times m}$ of the null space $N_1$ of
\begin{equation}
 \left[
                \begin{array}{c}
                 L_A(u)-L_A(0)\\
                 L_B(v)-L_B(0)
                 \end{array}\right]=0.
\end{equation}
The solutions $(\alpha_i,\beta_j)$ of~(\ref{SysStart}) arising
from~(\ref{alphas}) and~(\ref{betas}) correspond to $N_1\cap
(A\times B)$.

Fixing a particular solution $w_2$ of
\begin{equation}
  \bfA w+\bfB \bfP_{h^*}(\bfC w + \bfd)=0,
\end{equation}
choose a basis $W_2\in\bC^{2k\times m}$ of the null space $N_2$ of
\begin{equation}
  \bfA w+\bfB \bfP_{h^*} \bfC w=0.
\end{equation}
We have the intrinsic homotopy with variable $y\in\bC^m$
\begin{equation}\label{dumbHomotopy}
  \sF\left((1-\tau)\left[w_1+W_1y\right] + \tau
     \left[w_2+W_2y\right]\right)=0.
\end{equation}
Since $N_1$ is transverse to $A\times B$, the
$(2k-a-b)$-dimensional affine subspace given by
\begin{equation}
   \left\{
   \tau_1 \left[w_1+W_1y\right]
 + \tau_2 \left[w_2+W_2y\right]\ \big|\ y\in\bC^m
   \right\}
\end{equation}
is transverse to $A\times B$ for all but a finite set of
$[\tau_1,\tau_2]\in \pn 1$. In particular for all but a finite number of
$\gamma\in \bC$ of absolute value one,
with the relation between $\tau$ and $t$ as in~(\ref{crazyFormula}),
the $m$-dimensional affine subspace given by
\begin{equation}\label{yaTLS}
 \left\{
   (1-\tau)\left[w_1+W_1y\right] + \tau \left[w_2+W_2y\right]
   \ \big|\ y\in\bC^m
 \right\}
\end{equation}
is transverse to $A \times B$ for all $t\in(0,1]$. By genericity in
the choices of $\bfA,\bfB,\bfC,\bfd$, this is true for $t=0$ also.
Thus using the homotopy~(\ref{dumbHomotopy}) to track the
paths starting with the $(\alpha_i,\beta_j)$ at $t=1$, we get the
start solutions of the cascade at $t=0$.

In practice it will be convenient to go directly from solutions
$(\alpha_i,\beta_j)$ of~(\ref{SysStart}) arising
from~(\ref{alphas}) and~(\ref{betas}) to $\sE_{\hmax-1}$ or any
$\sE_j$  with $j<\hmax$.  Doing this we want to know that the
limits of the paths of the intrinsic homotopy starting with the
solutions $(\alpha_i,\beta_j)$ contain the subset
$\sG_j$ for which $z_j\not=0$ and a set of
points $\hatW_j$ which contains a set of
witness points $W_j$.  This is true for both the intrinsic and the
earlier extrinsic homotopy of \cite{SVW6}.  The reason why this is
so is that the solutions $\sG_j\cup\hatW_j$ are contained in the
set of isolated solutions of $\sE_j$ restricted to $A\times B$.
Therefore by \cite[Lemma A.1]{SVW6}, there is a Zariski open set
of $t\in\bC$ such that except for a finite choice of $\gamma $ of
absolute value one in (\ref{crazyFormula}), $\sG_j\cup\hatW_j$ are
limits of isolated solutions of the homotopy (\ref{dumbHomotopy})
restricted to $A\times B$. Since the solutions at $t=1$ of the
homotopy (\ref{dumbHomotopy}) on $A\times B$ are the transversal
intersection with the $m$-dimensional affine subspace given by
Eq.(\ref{yaTLS}), it follows that for the $t$ near $1$ this is
still true. Thus the isolated solutions of the homotopy
(\ref{dumbHomotopy}) for a Zariski open set of the $t$ are
continuations from solutions $(\alpha_i,\beta_j)$
of~(\ref{SysStart}) arising from~(\ref{alphas}) and~(\ref{betas}),
and in consequence $\sG_j\cup\hatW_j$ are contained in limits of
isolated solutions of the homotopy (\ref{dumbHomotopy}) restricted
to $A\times B$ starting at these points.

The current default is to go directly from solutions
$(\alpha_i,\beta_j)$ of~(\ref{SysStart}) arising
from~(\ref{alphas}) and~(\ref{betas}) to $\sE_{\hmax-1}$.

\subsection{Intrinsic Cascade Homotopies}\label{Sec:IntrinsicCascade}

In this section, we convert the extrinsic cascade homotopies
of~(\ref{Eq:CascadeHomotopies}) into intrinsic the form
of~(\ref{Eq:LinearIntrinsic}).  This must be done a bit more
delicately than what was done for the start homotopy, because we
must preserve the containment of $\sH_{i,j}$ inside the parameter
space of $\sE_i$ so that we retain the properties stated at the
end of \S\ref{Sec:ExtrinsicReview}.  We do this by deriving an
intrinsic homotopy whose path is exactly the same as a generic
real path  from $\tau=1$ to $\tau=0$ in~(\ref{Eq:CascadeHomotopies}).

We start by eliminating $z$ by substitution from the last block
row of~(\ref{Eq:CascadeHomotopies}) into the middle row.  We use
the facts that for $i>j$, $\bfP_i\bfP_j=\bfP_j$ and
$\bfP_i\bfP_{ji}=\bfP_{ji}$ to obtain
\begin{equation}\label{Eq:CascadeHomNoZ}
 \sH_{i,j}(t,w):=
 \left[
    \begin{array}{l}
       \sF(w)\\
       \bfA w + \bfB (\bfP_j + \tau \bfP_{ji})(\bfC w + \bfd)
    \end{array}
  \right]=0,
\end{equation}
which, abusing notation, we still call $\sH_{i,j}$. By similar
abuse of notation, we use $\sE_h(w)$ in place of $\sE_h(w,z)$
after eliminating $z$ from~(\ref{Eq:Cascade}).

Our first observation concerns the existence of a constant
particular solution throughout the cascade.
\begin{lemma}\label{Lem:Particular} The inhomogeneous
linear system
\begin{equation}
  \left[
    \begin{array}{c}
      \begin{array}{cc} I_k & -I_k \end{array} \\
      \bfC
    \end{array}
  \right] w =
  \left[
    \begin{array}{c}
       0 \\
       -\bfd
    \end{array}
  \right]
\end{equation}
has a unique nonzero solution $\epsilon$.
\end{lemma}

\begin{proof}The genericity of $\bfC$ implies the
invertibility of
$\displaystyle
  \left[
    \begin{array}{c}
      \begin{array}{cc} I_k & -I_k \end{array} \\
      \bfC
    \end{array}
  \right]$.
\end{proof}

Notice that this implies that both $\bfA \epsilon=0$ and $\bfC
\epsilon + \bfd=0$, and therefore $w=\epsilon$ is a solution of
\begin{equation}
  \bfA w + \bfB (\bfP_j + \tau \bfP_{ji})(\bfC w + \bfd)=0
\end{equation}
for any $i,j,\tau$.

Let $\bfY_h$ be the homogeneous linear system
\begin{equation} \label{Eq:defYh}
  \bfY_h:=(\bfA + \bfB \bfP_h \bfC) w = 0.
\end{equation}
The following lemma concerning the null space of $\bfY_h$ is
crucial for the conversion to an intrinsic form.
\begin{lemma}\label{Lem:NullSpaceY}
For any $j$ and $i$ such that $\hmin\le j<i\le \hmax$, there exist
matrices $E\in\bC^{2k\times(m-i+j)}$ and
$F,G\in\bC^{2k\times(i-j)}$ such that
\begin{enumerate}
  \item $[E\>\>F] = \Null \bfY_i$
  \item $[E\>\>G] = \Null \bfY_j$
  \item $\bfP_{ji}\bfC F = \bfP_{ji}\bfC G =
  \left[
   \begin{array}{c}
    0 \\
    I_{i-j} \\
    0
  \end{array} \right]$,
\end{enumerate}
where the $(i-j)\times(i-j)$ identity matrix $I_{i-j}$ appears in
rows $j+1,\ldots,i$.
\end{lemma}

\begin{proof}
We must first establish that $\bfY_i$ and $\bfY_j$ are full row
rank $a+b$ so that $m=2k-a-b$ is the correct dimension of their
null spaces.  Since $\bfA$ depends on generic $\bA$
(see~(\ref{Eq:bfA})) and $\bfB$ and $\bfC$ are generic, it suffices
to show that there is at least one choice of $\bA$, $\bfB$, $\bfC$
such that $\bfY_h$ is full rank for $\hmin\le h \le \hmax$.  For
$a+b<k$, it suffices to choose $\bfB=0$, $\bfC=0$ and choose $\bA$
to make $\bfY_h=[I_{a+b}\>\> 0\>\>-I_{a+b}\>\>0]w$.  For $a+b>k$,
choose $\bA=[I_k\>\>0]^T$, choose $\bfB$ with $I_{k-a-b}$ in the
lower left and $\bfC$ with $I_{k-a-b}$ in the upper left.  Since
$h\ge k-a-b$, this suffices to make $\bfY_h$ full rank, as one may
check by direct substitution.

Next, we establish that $\bfY_i$ and $\bfY_j$ share a null
subspace of dimension $m-i+j$.  Note that
\begin{equation}
  \bfY_i = (\bfA + \bfB(\bfP_j + \bfP_{ji})\bfC)w = (\bfY_j +
  \bfB\bfP_{ji}\bfC) w.
\end{equation}
The matrix $\bfB\bfP_{ji}\bfC$ is independent of
$\bfB\bfP_{j}\bfC$ because the projection matrices pick out
different rows and columns of generic matrices $\bfB$ and $\bfC$.
Accordingly, the subspace $\Null\bfY_i\cap\Null\bfY_j=\Null\bfY_j
\cap \Null (\bfB\bfP_{ji}\bfC)$.  These have dimension $m$ and
$(i-j)$, respectively, and they meet transversely, so the
intersection has dimension $m-i+j$. Let $E$ be any basis for this
subspace.

Now, suppose $\hat F$ completes a basis $[E\>\>\hat F]$ for
$\Null\bfY_i$. It must be independent of
$\Null(\bfB\bfP_{ji}\bfC)$, and since $\bfB$ is generic, this
implies that $\bfP_{ji}\bfC\hat F$ must be full rank.  Since
$\bfP_{ji}$ zeros out all but rows $j+1,\ldots,i$, this implies
that
\begin{equation}
  \bfP_{ji}\bfC \hat F = \left[ \begin{array}{c}
    0 \\
    Q \\
    0
  \end{array} \right]
\end{equation}
must have a full-rank $(i-j)\times(i-j)$ matrix $Q$ in rows
$j+1,\ldots,i$.  Then, $F=\hat F Q^{-1}$ completes the basis of
$\bfY_i$ while also satisfying Condition~3 of the lemma.  Similar
reasoning shows the existence of $G$.
\end{proof}

Choosing a random $\gamma\in \bC$, we form the linear system
\begin{equation}\label{Eq:W}
  W_{i,j}(t,y) = \epsilon +
           \left[\begin{array}{cc} E & tF+\gamma(1-t)G
           \end{array}\right]y
\end{equation}
where $y\in\bC^m$.  From this, we form the intrinsic homotopy
\begin{equation}\label{Eq:IntrinsicHomotopy}
   H_{i,j}(t,y) = \sF\left( W_{i,j}(t,y) \right) = 0,
\end{equation}
and track $y$ as $t$ goes from 1 to 0 on the real interval.

The crucial fact behind the equivalence of the intrinsic and
extrinsic homotopies is that the space intrinsically parameterized
in~(\ref{Eq:W}) is the same for appropriate choices of
parameters as the space that we extrinsically cut out with linear
equations before.
\begin{lemma}\label{lemma2.4Plus}
For all but a finite number of $\gamma\in \bC$ of absolute value
one, it follows that for any $t\in[0,1]$ there is a
$0\not=\tau\in\bC$ such that the kernel of the linear system
\begin{equation}\label{Eq:Equiv}
  \bfA w+ \bfB(\bfP_j +
  \tau\bfP_{ji})(\bfC w+\bfd) = 0.
\end{equation}
on $\bC^{2k}$ is parameterized by $W_{i,j}(t,y)$ where
$y\in\bC^m$.
\end{lemma}
\begin{proof} This follows immediately for $t=0$ and $1$ with no
restriction
on $\gamma$ of absolute value $1$ by taking $\tau$ equal to $0$
and $1$ respectively.

Combining this with the  dimension of the kernel
of~(\ref{Eq:Equiv}) being at least $m$, we conclude that the
dimension of the kernel of~(\ref{Eq:Equiv}) is exactly $m$
except for finitely many $0\not=\gamma\in\bC$.  In particular, for
all but a finite number $\gamma$ of absolute value $1$, the
dimension of the kernel of~(\ref{Eq:Equiv}) for
$\tau$ and $t$ as in~(\ref{crazyFormula})
with $t\in(0,1)$ is of dimension $m$. Since
$\epsilon$ satisfies both $\bfA\epsilon=0$ and
$\bfC\epsilon+\bfd=0$, it is therefore enough to show that for all
$(t,y)$
\begin{equation}
  \left( \bfA +\bfB(\bfP_j + \tau\bfP_{ji})\bfC \right)
  \left[\begin{array}{cc} E & tF+\gamma(1-t)G \end{array}\right]y=0.
\end{equation}

Since the columns of $E$ are in
$\Null\bfY_j\cap\Null(\bfB\bfP_{ji}\bfC)$, it is annihilated.
Since $y$ is arbitrary, we must have
\begin{equation}
  \left( \bfA +\bfB(\bfP_j + \tau\bfP_{ji})\bfC \right)
  [ tF+\gamma(1-t)G ] =0.
\end{equation}
Since $F$ is in $\Null\bfY_i$ and $G$ is in
$\Null\bfY_j$, this is the same as
\begin{equation}
  \bfB \left( (\tau-1)t\bfP_{ji}\bfC F
    + \tau\gamma(1-t)\bfP_{ji}\bfC G \right) =0.
\end{equation}
By Condition~3 of Lemma~\ref{Lem:NullSpaceY}, this becomes
\begin{equation}
  \left( (\tau-1)t + \tau\gamma(1-t) \right)
  \bfB
   \left[ \begin{array}{c}
    0 \\
    I_{i-j} \\
    0
  \end{array} \right],
\end{equation}
which equals zero by~(\ref{crazyFormula}).
\end{proof}

We rephrase Lemma \ref{lemma2.4Plus}.
\begin{lemma}\label{lemma2.4PlusVer2}
For all but a finite number of $\gamma\in \bC$ of absolute value
one, it follows that for any $t\in[0,1]$, the system
\begin{equation}
  \sF\left( W_{i,j}(t,y) \right) = 0
\end{equation}
on $\bC^m$ is the intrinsic system associated to the system
\begin{equation}
 \left[
    \begin{array}{l}
       \sF(w)\\
       \bfA w + \bfB \bfP_i z\\
       z - (\bfP_j + \tau \bfP_{ji})(\bfC w + L)
    \end{array}
  \right]=0
\end{equation}
with $\tau = t/(t+\gamma(1-t))$.
\end{lemma}

We define $\sG_i$ as the set of nonsingular solutions of
$H_{i,j}(1,\omega)$ on which $\bfP_i(\bfC w + L)$ is nonzero
and which correspond to points of $A\times B$;
$\sG_j$ as the set of nonsingular solutions of
$H_{i,j}(0,\omega)$ on which $\bfP_j(\bfC w + L)$ is nonzero
and which correspond to points of $A\times B$; and
$\sG_{i,j}$ as the sent of limits obtained by tracking $\sG_i$
from $t=1$ to $t=0$ using the homotopy $H_{i,j}(t,\omega)$.

\begin{theorem}\label{mainTheorem}The subset $\hatW_j\subset
\sG_{i,j}$ on which $\bfP_j(\bfC w + L)$ is zero contains a set of
witness points for the $j$-dimensional components of $A\cap B$.
These witness points include $\deg(Z)$ distinct points for each
irreducible $j$-dimensional component $Z$ of $A\cap B$. Moreover
$\sG_j\subset\sG_{i,j}$.
\end{theorem}

\begin{proof} The sets $\sG_i,\sG_j$ considered as sets of solutions
of the extrinsic systems $\sE_i,\sE_j$ on $\bC^{2k}$ are the same
as the sets occurring in the homotopy of \cite{SVW6}. The
extrinsic homotopy from \cite{SVW6} that we discussed in \S
\ref{Sec:ExtrinsicReview} is simply a differentiable path $P$
parameterized by $t\in[0,1]$ on a complex line $\ell$ in the
parameter space of the systems $\sE_i(w,z)$ joining a general
point $\sE_i$ to a general point $\sE_j$ of the linear subspace of
systems of the from $\sE_j(w,z)$.  The only fact about the path
$P$ used in \cite{SVW6} is that it depends on a choice of
$\gamma\in\bC$ of absolute value~$1$, which can be chosen, except
for a finite number of complex numbers of absolute value~$1$, so
that $P$ avoids a certain finite subset $B$ of $\ell$. In Lemma
\ref{lemma2.4PlusVer2}  we show that the intrinsic homotopy leads
to systems on the {\em same} complex line $\ell$.  What changed is
that the path $P'$ on $\ell$ is not linearly related to the
original path $P$.  But since the path $P'$ depends on a choice of
$\gamma\in\bC$ of absolute value $1$, which can still be chosen,
except for a finite number of complex numbers of absolute value
$1$, so that $P'$ avoids the  finite subset $B$ of $\ell$, the
same conclusions of \cite{SVW6} still hold.
\end{proof}

\subsection{Algorithm Summary}

The homotopy algorithm to intersect two positive dimensional
varieties in intrinsic coordinates is described below.
After the initialization, there are three stages.
First is the homotopy to start the cascade, followed
by the homotopy to find a witness sets for the top dimensional
part of $A \cap B$.  Thirdly, all lower dimensional parts of
$A \cap B$ are computed in a loop from $b-2$ down to~$h_0$.
The second and third stage are separate because we can avoid
a coordinate transformation.  Also, in many cases -- such as the
important application of the intersection with a hypersurface --
the loop will never be executed.

Some subroutines used in the algorithm below are just implementations
of one formula in the paper, e.g.: {\bf Combine} implements~(\ref{Eq:sFw}).
Next we describe briefly the other subroutines.

The linear algebra operations to deal with solutions in intrinsic
coordinates are provided in the subroutines {\bf Start\_Plane},
{\bf Project}, {\bf Initialize}, {\bf Basis}, and {\bf Transform}.
Given the equations for $L_A$ and $L_B$, {\bf Start\_Plane} first
computes a basis for the null space of $L_A^{-1}(0)$ and
$L_B^{-1}(0)$ before doubling the coordinates into a corresponding
basis in $\bC^{2k}$. After orthonormalization of the basis, {\bf
Project} computes the intrinsic coordinates for the product of the
given witness sets of $A$ and $B$.  The subroutine {\bf
Initialize} first generates the random matrices $\bfA$, $\bfB$,
$\bfC$, and $\bf d$ before computing the $\epsilon$ of
Lemma~\ref{Lem:Particular}. In addition, {\bf Initialize} returns
the operator $\bfY$, which returns for any $h$ the corresponding
$\bfY_h$ of~(\ref{Eq:defYh}). Lemma~\ref{Lem:NullSpaceY} is
implemented by {\bf Basis}, while {\bf Transform} converts the
coordinates for the solutions from one basis into another.

The path tracking is done by the procedure {\bf Track}.
On input are the homotopy and start solutions.  Except from the
set up of the homotopy in intrinsic coordinates, one can implement
{\bf Track} along the lines of general path following methods,
see~\cite{AG90}, \cite{Li97,Li03}, or~\cite{Mor87}.

The subroutine {\bf Filter} takes on input the witness sets $\sW$ for
higher dimensional components and the list $\sZ$.  On return is $\sW$,
augmented with a witness set for the solution set at the current
dimension, and a filtered list $\sZ$ of nonsolutions.
The list $\sZ$ given to {\bf Filter} may contain points on higher
dimensional solution sets.  To remove such points, a homotopy
membership test as proposed in~\cite{SVW2} can be applied.
Recently, an interesting alternative was proposed
by Li and Zeng in~\cite{LZ}.
The nonsolutions serve as start solutions in the cascade to find
witness sets for the lower dimensional solution sets.
If $\sZ$ becomes empty after {\bf Filter}, the algorithm terminates.

\begin{algorithm} {\rm
Intersecting two Positive Dimensional Varieties $A$ and $B$.

\noindent ~\hspace{-0.3cm}~\begin{tabular}{lr}
 Input: $k$, $a$, $b$ $(a \geq b)$;
   &  {\em $\dim(A) = a$, $\dim(B) = b$, $A,B \subset \bC^k$} \\
\hspace{1.1cm} $f_A(u) = 0, f_B(v) = 0$;
   &  {\em polynomial systems in $u,v \in \bC^k$} \\
\hspace{1.1cm} $L_A(u) = 0, L_B(v) = 0$;
   &  {\em $\dim(L_A^{-1}(0)) = k\!-\!a$, $\dim(L_B^{-1}(0)) = k\!-\!b$} \\
\hspace{1.1cm} $\sW_A, \sW_B$.
   & {\em solutions in witness sets for $A$ and $B$} \\
 Output: $\sF(x) = 0$;
   &  {\em system combined from $f_A$, $f_B$ in $x \in \bC^k$} \\
\hspace{1.4cm} $L = [L_{h_0}, \ldots , L_{b-1}]$;
   &  {\em list of linear spaces, $\dim(L_i^{-1}(0)) = i$} \\
\hspace{1.4cm} $\sW = [\sW_{h_0}, \ldots , \sW_{b-1}]$.
   &  {\em solutions $\sW_i$ in $i$-dim witness sets}
\end{tabular}

\noindent ~\hspace{-0.5cm}~\begin{tabular}{lr}
$\sF := {\bf Combine}(f_A,f_B)$;
   &  {\em combine systems $f_A$ and $f_B$ as in}~(\ref{Eq:sFw}) \\
$S := {\bf Start\_Plane}(L_A,L_B)$;
   &  {\em basis for plane defining $\sW_A \times \sW_B$} \\
$\sZ := {\bf Project}(\sW_A \times \sW_B,S)$;
   &  {\em solutions to start the cascade} \\
$[\bfY,\epsilon] := {\bf Initialize}(k,a,b)$;
   &  {\em linear space $\bfA w + \bfB \bfP_h \bfC(w +\bfd) = 0$} \\
$[E,F,G] := {\bf Basis}(\bfY_{b},\bfY_{b-1})$;
   &  {\em basis for $\Null \bfY_b$ and $\Null \bfY_{b-1}$} \\
$W(t,y) := [ t S + (1-t)[\epsilon + [E~F]]$;
   &  {\em deform start plane $S$ into $[E~F]$} \\
   &  {\em with $t$ using formula~}(\ref{crazyFormula}) \\
$\sZ := {\bf Track}(\sF(W(t,y)),\sZ)$;
   &  {\em homotopy to start the cascade} \\
$\sZ := {\bf Track}(\sF,[E,F,G],\sZ)$;
   &  {\em find top dimensional component} \\
$[\sW_{b-1},\sZ] := {\bf Filter}(\sW,\sZ)$;
   &  {\em keep witness sets and nonsolutions} \\
$h_0 := \max(a+b-k,0)$;
   &  {\em minimal $\dim(A\cap B)$} \\
for $j$ from $b-2$ down to $h_0$ do
   &  {\em compute witness set at dimension $j$} \\
\hspace{0.3cm} $[E,F,G] := {\bf Basis}(\bfY_{j+1},\bfY_j)$;
   &  {\em $W(t,y) = \epsilon
           + [E ~~ t F+\gamma(1-t)G] y$} \\
\hspace{0.3cm} $\sZ := {\bf Transform}(\sZ,[E,F])$;
   &  {\em coordinates into new basis~$[E~F]$} \\
\hspace{0.3cm} $\sZ := {\bf Track}(\sF,[E,F,G],\sZ)$;
   &  {\em homotopy $\sF(W_{j+1,j}(t,y)) = 0$} \\
\hspace{0.3cm} $[\sW_j,\sZ] := {\bf Filter}(\sW,\sZ)$;
   &  {\em keep witness sets and nonsolutions} \\
end for. \\
\end{tabular}
}
\end{algorithm}

\section{Numerical Experiments} \label{Sec:comparison}

The algorithms in this paper have been implemented and tested
with PHCpack~\cite{V99}.  To compare with our implementation in
extrinsic coordinates, we use the same examples as in~\cite{SVW6}.
All computations were done on a 2.4 Ghz Linux machine.

\begin{description}
\item[(1) An Example from Calculus.] In this example, we intersect
a cylinder~$A$ with a sphere~$B$.  More precisely,
$A = \{ \ (x,y,z) \ | \  x^2 + y^2 - 1 = 0 \ \}$
and $B = \{ (x,y,z) \ | \  (x+0.5)^2 + y^2 + z^2 - 1 = 0 \ \}$.
The intersection $A \cap B$ is a curve of degree four.
Since $k = 3$, $a = 2$, and $b = 2$: $h_0 = 1$, so there are
only two homotopies, each defining four solution paths.

\item[(2) An Illustration of the Cascade.]  In this example we
need to execute the cascade to find the point of intersection.
We consider the components $A = \{\ x = 0, y=0 \ \}$ and
$B = \{ \ z = 0, w =0 \ \}$ as solution sets of the same
system $f(x,y,z,w) = [ xz , xw, yz, yw ]^T = 0$.
We have $k=4$, $a=2$, and $b=2$.

\item[(3) Adding an Extra Leg to a Moving Platform.]  In this
example we cut a hypersurface $A$ in $\bC^8$ with a curve $B$,
i.e.: $a = 7$ and $b = 1$.  The application concerns a
Griffis-Duffy platform~\cite{GD93} (analyzed by Husty and
Karger in~\cite{HK00} and subsequently in~\cite{SVW5})
where $A \cap B$ can be interpreted as adding a seventh
leg to the platform so it no longer moves.
As $\deg(A) = 2$ and $\deg(B) = 28$ (ignoring the mechanically
irrelevant components), there are 56 paths to trace,
by two homotopies.

\end{description}
In the Table~\ref{tabdim} below we list all important dimensions of the
three example applications.  A summary of the execution times is
reported in Table~\ref{tabtimings}.

\begin{table}[hbt]
\begin{center}
\begin{tabular}{|c||c||c|c||c|c||c|c|c|} \hline
  \multicolumn{2}{|c||}{example}
  & \multicolumn{4}{c||}{dimensions and degrees of $A$ and $B$}
  &  & & $\deg(A)$~~~ \\ \cline{2-6}
     & $k$ & $\dim(A)$ & $\deg(A)$ & $\dim(B)$ & $\deg(B)$
     & $m$ & $M$ & $\times \deg(B)$ \\ \hline
 (1) & 3 & 2 & 2 & 2 & ~2 & 2 & ~7 & ~4 \\
 (2) & 4 & 2 & 1 & 2 & ~1 & 4 & 10 & ~1 \\
 (3) & 8 & 7 & 2 & 1 & 28 & 8 & 17 & 56 \\ \hline
\end{tabular}
\caption{\label{tabdim}
         Dimension and degrees of the two irreducible sets $A$ and $B$
         for the three examples, followed by \#variables
         $m=2k-\dim(A)-\dim(B)$, $M = 3k-a$ (which is the \#variables in
         the extrinsic homotopy), and number of paths $\deg(A) \times \deg(B)$
     at the start of the cascade.}
\end{center}
\end{table}

\begin{table}[hbt]

\begin{center}
\begin{tabular}{|c||c|c|c||c|c|} \hline
     & \multicolumn{3}{c||}{Homotopies}
     & \multicolumn{2}{c|}{Total CPU Time} \\ \cline{2-6}
     &   0  &   1  &   2  & intrinsic & extrinsic \\ \hline \hline
 (1) & 0.03 & 0.01 &  --  & ~0.04 & ~0.07  \\
 (2) & 0.01 & 0.02 & 0.01 & ~0.04 & ~0.11 \\
 (3) & 9.90 & 5.94 &  --  & 15.84 & 34.70 \\ \hline
\end{tabular}
\caption{\label{tabtimings}
        Timings in CPU user seconds on 2.4Ghz Linux machine.
         The second column concerns the homotopy to start the cascade,
         in the third column are the timings for the top dimensional
         components, followed by the eventual next homotopy in the
         cascade.}
\end{center}
\end{table}

In these numerical experiments, we save about half of the
computational time when working in intrinsic coordinates.
Comparing the number of variables of the original extrinsic
method, $M = 3k - a$ for the examples tested,
with the number for the intrinsic method,
$m=2k-\deg(A)-\deg(B)$, we have in these experiments $3k-a=7,10,17$
variables reduced to $2,4,8$, or more than half.  Since the cost of
linear solving is ${\cal O}(n^3)$, this implies about a eight-fold
reduction in the cost of that portion of the code.  Linear solving
can be a significant portion of the total cost, as it is used in
Newton's method for tracking the homotopy paths.  The experimental
results suggest that this was accounting for about half of the
total cost in the extrinsic method, but accounts for a much less
significant fraction of the computational cost of the intrinsic
method.  The other 50\% or so of the cost remains, which is
attributable to function evaluation, data transfer, and other
overhead.  The cost of function evaluation can vary dramatically
from one polynomial system to another, so we cannot definitively
expect the same percentage savings for all systems, but we can say
that the intrinsic formulation seems to give a substantial
reduction in computational time.

\end{document}